\documentstyle[12pt]{article}
\begin{document}
\begin{center}
{\bf UNIFICATION OF INDEPENDENCE \\
IN QUANTUM PROBABILITY}\\[40pt]
{\sc Romuald Lenczewski} \\[40pt]
{\it Hugo Steinhaus Center for Stochastic Methods\\
Institute of Mathematics, Technical University of Wroc{\l}aw\\
50-370 Wroc{\l}aw, Poland}\\
e-mail lenczew@im.pwr.wroc.pl\\[20pt]
\end{center}
\begin{abstract}
Let ($*_{l\in I}{\cal A},*_{l\in I}$($\phi_{l}, \psi_{l}$)),
be the conditionally free product of unital free *-algebras ${\cal
A}_{l}$, where $\phi_{l}, \psi_{l}$ are states on ${\cal A}_{l}$, 
$l\in I$.
We construct a sequence of
noncommutative probability spaces ($\widetilde{{\cal A}}^{(m)}, \widetilde{\Phi}^{(m)}$), 
$m\in {\bf N}$, where
$\widetilde{{\cal A}}^{(m)}=\bigotimes_{l\in I}
\widetilde{{\cal A}}_{l}^{\otimes m}$ and
$\widetilde{\Phi}^{(m)}=\bigotimes_{l\in I}\widetilde{\phi}_{l}\otimes 
\widetilde{\psi}_{l}^{\otimes (m-1)}$, $m\in {\bf N}$, 
$\widetilde{{\cal A}}_{l}={\cal A}*{\bf C}[t]$,  
and the states $\widetilde{\phi}_{l}, \widetilde{\psi}_{l}$ are Boolean
extensions of $\phi_{l}, \psi_{l}$, $l\in I$,
respectively. We define unital *-homomorphisms $j^{(m)}$
$:*_{l\in I}{\cal A}_{l}\rightarrow \widetilde{{\cal A}}^{(m)}$ such that
$\widetilde{\Phi}^{(m)}\circ j^{(m)}$ converges pointwise to
$*_{l\in I}(\phi_{l},\psi_{l})$. Thus, the variables
$j^{(m)}(w)$, where $w$ is a word in $*_{l\in I}{\cal A}_{l}$,   
converge in law to the conditionally free variables. 
The sequence of noncommutative probability spaces
(${\cal A}^{(m)},\Phi^{(m)}$), where ${\cal
A}^{(m)}=j^{(m)}(*_{l\in I}{\cal A}_{l})$ and
$\Phi^{(m)}$ is the restriction of $\widetilde{\Phi}^{(m)}$ to
${\cal A}^{(m)}$, is called a {\it hierarchy of freeness}. Since all
finite joint correlations for known examples of independence
can be obtained from tensor products of appropriate
*-algebras, this approach can be viewed as a unification of independence.
Finally, we show how to make the $m$-fold free product 
$\widetilde{{\cal A}}^{*(m)}$ into
a cocommutative *-bialgebra associated with $m$-freeness.\\[10pt] 
\end{abstract}
\begin{center}
{\sc 1. Introduction}
\end{center}
The aim of this paper is to show that the main types of noncommutative
independence can be
obtained from tensor independence and are related to 
appropriately constructed *-bialgebras. 

Since different models have led to
almost separate theories and techniques, it seems desirable to
develop one theory covering all the cases, including
tensor, free and Boolean independence as well as their various
modifications. This work makes the first step in this
direction, namely is provides a unified treatment of the main
notions of independence existing in the literature in the sense
that it reduces the problem of calculating finite joint
correlations to a similar problem formulated for the tensor
product of *-algebras.
In other words, it shows that the main types of products of states
can be reduced to tensor products of states or they are pointwise 
limits of such states like in the case of freeness.

From the axiomatic approach presented in [Sch2] it follows that
under certain assumptions there are three 
``pure'' kinds of independence, namely commutative [C-H, G-vW], Boolean [vW]
and free [Voi, V-D-N], each with characteristic combinatorics. An
interpolation between the Boolean model and the free 
model has been given in terms of the so-called conditional
freeness (earlier called $\psi$-independence) [B-S], which allows
us to extract both models as special cases. Essentially, the
conditionally free probability is based on the approach to the free
probability presented by Voiculescu. Noncommutative probability
spaces obtained from this kind of approach have always been
viewed as very noncommutative and thus not directly related to
the tensor product case. Our theory provides a sequence of explicit
tensor product constructions which allows to (pointwise) approximate
the conditionally free product of states and thus may be viewed
as a unifying tensorization scheme. 

In our approach, instead of making the theory noncommutative in the
definition of the product of *-algebras, we stick to the
tensor product and simply take noncommutative extensions of
those *-algebras with non-canonical embeddings. 
Our ideas go back to the central limit theorem for the *-Hopf algebra
$U_{q}(su(2))$ in [Len1, Len2]. This was our version of
``$q$-independence'', which gave the $q$-Gaussian law in the
$q$-central limit theorem. It seemed natural that one should be
able to construct suitable *-bialgebras
associated with free independence and Boolean independence. An axiomatic
approach to this subject was
presented by Sch\"{u}rmann [Sch2].
Our approach is different and is the first one which gives explicit
*-homomorphic embeddings of the free product of unital free *-algebras 
into suitable tensor products. 
Moreover, our construction can also be used for
*-algebras for which ${\cal A}^{0}$ is a *-subalgebra of ${\cal A}$, where
${\cal A}={\cal A}^{0}\oplus {\bf C}{\bf 1}$. Thus, it is not less
general than the approach in [Sch2] (see Section 3).
Moreover, it gives a nice structure embodied by
the constructed hierarchy of freeness and this way fills
the ``gap'' between Boolean independence and freeness.

The main idea consists in constructing *-bialgebras (or,
*-Hopf algebras, if possible) similar to the $q$-bialgebras 
or the $q$-deformed enveloping algebras $U_{q}(su(2))$, but perhaps
``more noncommutative'', i.e. with ``more noncommutative''
kernels replacing the kernels given by $q$-relations studied in
[Sch1, Len1, Len2]. It presents no difficulty to construct a *-bialgebra
associated with Boolean independence, but in order to cover
freeness as well as the general case of conditional freeness,
one needs to construct a sequence of *-bialgebras in order to
obtain freeness as the limit in law (by which we understand the
convergence of finite joint correlations). 

This works for independent copies of the same
algebra. 
If we want to consider (free [Av, Voi], conditionally free [B-L-S],
Boolean [vW]) products of different algebras, a natural
generalization of the *-bialgebra techniques can be used.
Instead of the sequence of coproducts $\Delta^{(m)}$,
we take a sequence of *-homomorphisms $j^{(m)}$
(see Definition 2.1). 
It turns out that in order to obtain the Boolean product it is
enough to consider the 1-fold tensor product. The 2-fold tensor
product construction gives a noncommutative probability space that we
associate with 2-freeness, and so on, the $m$-fold tensor product giving
$m$-freeness. Consequently, in the limit $m\rightarrow
\infty$ we obtain freeness. In fact, all those
constructions can be embedded into one, using the infinite
tensor product of *-algebras, but it is convenient in some
places to carry out the proofs for the sequence of $m$-fold
tensor product constructions.

The implications of this fact should
lead to some new interesting developments of the theory. It
is not clear at this point to what extent our result will 
facilitate a unified approach to other aspects of quantum
probability.  
It is also hard to claim that the main results in quantum
probability will be reducible to the tensor product techniques
and, in the case of independent copies of the same *-algebra, to
the probability theory for *-bialgebras or *-Hopf algebras, no matter
how nice this connection might seem. However, we think that our result
provides a nice structure of independence 
in the noncommutative probability theory and perhaps will
lead to a unified treatment of such topics as limit theorems,
invariance principles, Fock spaces, etc.

In Section 2 we give basic definitions related to the extensions
of states on unital free *-algebras and we introduce a sequence
($\widetilde{{\cal A}}^{(m)}, \widetilde{\Phi}^{(m)}$) of quantum
probability spaces.
Namely, for each $m\in {\bf N}$ and given two unital free *-algebras 
${\cal A}_{1}$ and ${\cal A}_{2}$ we define the algebraic
tensor product 
$$
\widetilde{{\cal A}}^{(m)}=\widetilde{{\cal A}}_{1}^{\otimes m} \otimes
\widetilde{{\cal A}}_{2}^{\otimes m} 
$$
where $\widetilde{{\cal A}}_{l}={\cal A}_{l}*{\bf C}[t]$, $l=1,2$,
with hermitian $t$.
Given two pairs of states on ${\cal A}_{1}, {\cal A}_{2}$, namely
($\phi_{1},\psi_{1}$) and ($\phi_{2}, \psi_{2}$),
respectively, we construct the tensor product state 
$$
\widetilde{\Phi}^{(m)}\equiv \widetilde{\Phi}_{1}^{(m)}\otimes \widetilde{\Phi}_{2}^{(m)}
=\left(
\widetilde{\phi}_{1}\otimes
\widetilde{\psi}_{1}^{\otimes (m-1)}
\right)\otimes
\left(
\widetilde{\phi}_{2}\otimes \widetilde{\psi}_{2}^{\otimes (m-1)}
\right), 
$$
where $\widetilde{\phi}_{1}$, $\widetilde{\phi}_{2}$ are Boolean 
extensions of $\phi_{1}$, $\phi_{2}$, respectively (see
Definition 2.0), to states on
$\widetilde{\cal A}_{l}$, $l=1,2$.
For each $m\in {\bf N}$ we construct a *-subalgebra ${\cal
A}^{(m)}$ of $\widetilde{{\cal A}}^{(m)}$ on which 
the restriction of $\widetilde{\Phi}^{(m)}$ denoted by $\Phi^{(m)}$
can be interpreted as the (conditionally) $m$-free product state.
The pair (${\cal A}^{(m)}, \Phi^{(m)}$) is then
the noncommutative probability space associated
with $m$-freeness. In particular, 
$1$-freeness is in this scheme assigned to the Boolean product and
Boolean independence. 

In Section 3 we prove a number of technical results, especially
certain factorization lemmas which enable us to formulate our
main results.

These are presented in Section 4, where
we show that $\Phi^{(m)}\circ j^{(m)}$  
converges pointwise to the conditionally free product of states
$*_{l\in \{1,2\}}(\phi_{i},\psi_{i})$ on $*_{i\in \{1,2\}}{\cal A}_{i}$.
In particular, when $\psi_{l}=\phi_{l}$, $l=1,2$, we obtain in
the limit the free product of Voiculescu.
We also show how our results can be extended to the
case of infinitely many *-algebras.
An uncountable number of free *-algebras can be treated along
the same lines. 

In Section 5 we restrict ourselves to the case of one unital free *-algebra:
${\cal A}_{l}={\cal A}$ for all $l\in {\bf N}$. 
This corresponds to the case of conditionally free convolution
powers of states on ${\cal A}$. For each $m\in {\bf N}$ we equip 
the $m$-fold free product 
$\widetilde{{\cal A}}^{*(m)}=\widetilde{{\cal A}}* \ldots *\widetilde{{\cal
A}}$ ($m$ times)
with a *-bialgebra structure 
$(\widetilde{{\cal A}}^{*(m)}, \Delta^{(m)},\epsilon^{(m)})$ 
with coproduct 
$\Delta^{(m)}$ and counit $\epsilon^{(m)}$ 
(in this notation the symbols of products and units are supressed), 
which has an interesting property. Namely, if we lift two tensor
product states 
$\widetilde{\Phi}_{1}^{(m)},\widetilde{\Phi}_{2}^{(m)}$ to
states $\widehat{\Phi}_{1}$, $\widehat{\Phi}_{2}$, respectively,
on $\widetilde{{\cal A}}^{*(m)}$,
the convolution of $\widehat{\Phi}_{1}^{(m)}$ and
$\widehat{\Phi}_{2}^{(m)}$, which by definition is expressed
in terms of the coproduct as
$
\widehat{\Phi}_{1}^{(m)} \star \widehat{\Phi}_{2}^{(m)} \equiv
(\widehat{\Phi}_{1}^{(m)}\otimes \widehat{\Phi}_{2}^{(m)})\circ
\Delta^{(m)}, 
$
satisfies
$$
\lim_{m\rightarrow\infty}
\left(\widehat{\Phi}_{1}^{(m)}\star \widehat{\Phi}_{2}^{(m)}\right) 
\circ \widehat{i}_{1}(w) = 
(\phi_{1},\psi_{1})\star (\phi_{2}, \psi_{2})(w)
$$
where $\widehat{i}_{1}$ 
is the canonical *-homomorphic embedding of ${\cal A}$ into 
$\widetilde{{\cal A}}^{*(m)}$ given by 
$a\rightarrow a_{(1)}$, where $a_{(1)}$ is
the first copy of the generator $a$ in $\widetilde{{\cal
A}}^{*(m)}$. 

We view our model as a unified model of independence in the sense that
finite joint correlations for
known types of independence can be obtained from tensor
products of appropriately defined *-algebras and tensor
products of states. Here, the model of free probability of
Voiculescu takes the distinguished place of a limit case.
In a subsequent paper we will show that using infinite tensor
products and the GNS construction one can in fact embed all
levels of freeness in one tensor product of algebras.
A connection with the free product representation will alos be
established there.
\\[5pt] 
\begin{center}
{\sc 2. Preliminaries}
\end{center}
By a noncommutative probability space we understand a pair
$({\cal A}, \phi)$, where ${\cal A}$ is a unital *-algebra and
$\phi :{\cal A}\rightarrow {\bf C}$ is a state, i.e. a
normalized ($\phi({\bf 1})=1$), positive ($\phi(xx^{*}) \geq 0$ for
all $x\in {\cal A}$) functional.

Our construction will be carried out for unital free *-algebras ${\cal
A}$ generated by a set ${\cal G}^{+}$. We denote ${\cal G}^{-}=
\{a^{*}| a\in {\cal G}^{+}\}$, ${\cal G}={\cal G}^{+}\cup{\cal
G}^{-}$. Nonempty words in ${\cal A}$ will be denoted by $w=a_{1}\ldots
a_{k}$, where $a_{i}\in {\cal G}$. The length of $w$ will be
denoted by $l(w)$. We allow the empty word, which is denoted by
${\bf 1}$, of length $l({\bf 1})=0$.
The involution is given by the antilinear 
extension of $(a_{1}\ldots a_{k})^{*}=a_{k}^{*}\ldots a_{1}^{*}$.

For a given unital free *-algebra ${\cal A}$ we consider the free
product of ${\cal A}$ and ${\bf C}[t]$, the algebra of
polynomials in one hermitian variable $t$, which we denote 
$$
\widetilde{{\cal A}}= {\cal A}*{\bf C}[t].
$$
In this free product we identify units.
Also, we equip $\widetilde{\cal A}$ with a natural
involution defined by the antilinear extension of
$$
(t_{0}w_{1}t_{1}\ldots w_{n}t_{n})^{*}=t_{n}w_{n}^{*}\ldots
t_{1}w_{1}^{*}t_{0}, 
$$
where $w_{1}, \ldots, w_{n}$ 
are non-empty words in ${\cal A}$, and $t_{0}, \ldots , t_{n}$
are monomials in ${\bf C}[t]$,
respectively, of which $t_{1}, \ldots , t_{n-1}\neq {\bf 1}$.

Below we will define an extension of a state $\phi$ on
${\cal A}$  
to a state $\widetilde{\phi}$ on $\widetilde{{\cal A}}$ which we refer
to as the {\it Boolean extension} of $\phi$.\\
\indent{\par}
{\sc Definition 2.0.}
{\it For a given state ${\phi}$ on ${\cal A}$, we define a
Boolean extension of $\phi$ to be a functional 
$\widetilde{\phi}$ on $\widetilde{\cal A}$, which is the linear extension of
$\widetilde{\phi}({\bf 1})=1$ and} 
$$
\widetilde{\phi}\left(
t_{0}w_{1}t_{1}\ldots w_{n}t_{n}
\right)
= \phi(w_{1})\ldots\phi(w_{n})
$$
{\it where $w_{1}, \ldots, w_{n}$ are non-empty words in ${\cal A}$
and $t_{0}, \ldots, t_{n}$ are words in ${\bf C}[t]$,
of which $t_{1},\ldots , t_{n-1}$ are non-empty.}

One can obtain $\widetilde{\phi}$ from the Boolean product of $\phi$
and a *-homomorphism $h:{\bf C}[t]\rightarrow {\bf C}$, for
which $h(t)=1$. In fact, from the definition of the Boolean
product $\phi *_{B} h$ (see, for instance [B-L-S]), we obtain
$$
{\phi}*_{B}h\left(
t_{0}w_{1}t_{1}\ldots w_{n}t_{n}
\right)
= h(t_{0})\ldots h(t_{n})\phi(w_{1})\ldots\phi(w_{n})
$$
and using the assumptions on $h$ given above, we
obtain the Boolean extension of $\phi$.

From [B-L-S] it follows that $\widetilde{\phi}$ is a state.
It is also easy to see that the two sided *-ideal
generated by $t({\bf 1}-t)$ is contained in ${\rm
ker}\widetilde{\phi}$. Thus we can put 
$t^{n}=t$ and ${\bf 1}-t=({\bf 1}-t)^{n}$ in all formulas written modulo 
${\rm ker}\;\widetilde{\phi}$.
In other words, $\widetilde{\phi}$ does not distinguish between
positive powers of $t$.

One can say that the generator $t$ serves as a
``Boolean identity'', in contrast to 
$U_{q}(su(2))$-type Hopf algebras, where a similar object
satisfies certain $q$-commutation relation and can be viewed as
a ``$q$-identity''. Note that it plays
the role of a ``separator'' of words from the *-algebra ${\cal
A}$. This nice property will be crucial in further
considerations. 

Let us also recall the definition of the conditionally free product of
*-algebras. 
For a given family of unital *-algebras ${\cal A}_{l}$, $l\in
I$, and given pairs of states $\phi_{l}, \psi_{l}$ on ${\cal
A}_{l}$, one 
can define a state $\phi=*_{l\in I}(\phi_{l},\psi_{l})$ on their
free product $*_{l\in I}{\cal A}_{l}$ by $\phi({\bf 1})=1$ and 
the factorization property
$$
\phi(a_{1}\ldots a_{n})=\phi_{k_{1}}(a_{1})\ldots \phi_{k_{n}}(a_{n}),
$$
whenever $a_{j}\in {\cal A}_{k_{j}}$ and
$\psi_{k_{j}}(a_{j})=0$, where $k_{1}\neq k_{2} \neq \ldots \neq
k_{n}$. 
In particular, when $\psi_{j}=\phi_{j}$, we obtain the free
independence, and when $\psi_{j}=\pi_{1}$, where $\pi_{1}({\bf 1})=1$
and $\pi_{1}(w)=0$ for any non-empty word $w$, we get
Boolean independence. 

For given two unital free *-algebras ${\cal A}_{1}, {\cal A}_{2}$
generated by ${\cal G}_{1}^{+}, {\cal G}_{2}^{+}$, respectively,
let ${\cal G}_{l}={\cal G}_{l}^{+}\cup {\cal G}_{l}^{-}$,
where ${\cal G}_{l}^{-}=\{a^{*}|a\in {\cal G}_{l}^{+}\}$,
$l=1,2$. Given two pairs of states on those *-algebras, namely
$(\phi_{l},\psi_{l})$, $l=1,2$, we construct their
Boolean extensions 
$(\widetilde{\phi}_{l}, \widetilde{\psi}_{l})$ on $\widetilde{\cal A}_{l}$
as explained in Section 2. Using them, we will
construct for each $m \in {\bf N}$ a new noncommutative probability space
$(\widetilde{{\cal A}}^{(m)}, \widetilde{\Phi}^{(m)})$, where 
$$
\widetilde{{\cal A}}^{(m)}=\widetilde{{\cal A}}_{1}^{\otimes m}\otimes
\widetilde{{\cal A}}_{2}^{\otimes m}
$$
and the state $\widetilde{{\Phi}}^{(m)}$ is given by
$$
\widetilde{\Phi}^{(m)}=\widetilde{\phi}_{1}\otimes \widetilde{\psi}_{1}^{\otimes
(m-1)} \otimes
\widetilde{\phi}_{2}\otimes \widetilde{\psi}_{2}^{\otimes (m-1)}.
$$
The involution on the $2m$-fold tensor product is given by
$$
(b_{1}\otimes \ldots  \otimes b_{m}
\otimes c_{1}\otimes \ldots \otimes c_{m})^{*}
=b_{1}^{*}\otimes \ldots \otimes b_{m}^{*} 
\otimes c_{1}^{*} \otimes \ldots \otimes c_{m}^{*}.
$$
Let $i_{k,m}$, $m\in {\bf N}, k\in [m]\equiv\{1, \ldots , m\}$ 
be the canonical
*-homomorphic embeddings of $\widetilde{{\cal A}}_{l}$ into
$\widetilde{{\cal A}}_{l}^{(m)}$ (for each $l$ we use the same
notation), i.e.
$$
i_{k,m}(a)=I_{k-1}\otimes a \otimes I_{m-k},\;\;
i_{k,m}(t)=I_{k-1}\otimes t \otimes I_{m-k},
$$
where $a\in {\cal G}$, $I_{k}={\bf 1}^{\otimes k}$,
extended by linearity and multiplicativity to $\widetilde{{\cal
A}}_{l}$. We will adopt the convention that $i_{m+1,m}(a)=0$.
We will also use the abbreviated notation for products of
$i_{k,m}(t)$'s. Namely
$$
t_{[k,m]}=i_{k,m}(t)\ldots i_{m,m}(t)=I_{k-1}\otimes
t^{\otimes (m-k+1)}.
$$
\indent{\par}
{\sc Definition 2.1}
{\it For given $a\in {\cal G}_{1}, b\in {\cal G}_{2}$ and $m\in
{\bf N}$ let}
$$
j_{1}^{(m)}(a)=\sum_{k=1}^{m}(i_{k,m}(a)-i_{k+1,m}(a))\otimes
t_{[k,m]} 
$$
$$
j_{2}^{(m)}(b)=\sum_{k=1}^{m}t_{[k,m]}\otimes
(i_{k,m}(b)-i_{k+1,m}(b)) 
$$
{\it and define the *-homomorphism}
$$
j^{(m)}: {\cal A}_{1}*{\cal A}_{2}\rightarrow 
\widetilde{{\cal A}}^{(m)}
$$
{\it as the linear extension of $j^{(m)}({\bf 1})=I_{m}\otimes I_{m}$ and}
$$
j^{(m)}(w_{1}\ldots w_{n})=j_{k_{1}}^{(m)}(w_{1})\ldots 
j_{k_{n}}^{(m)}(w_{n}),
$$
{\it where $w_{1}, \ldots , w_{n}$ are non-empty words in
${\cal A}_{k_{1}}, \ldots , {\cal A}_{k_{n}}$,
where $k_{1},\ldots, k_{n}\in \{1,2\}$.}

Equivalently, we can write the above
condition in terms of the generators, i.e. 
$$
j^{(m)}(a_{1}\ldots a_{n})=j_{k_{1}}^{(m)}(a_{1})\ldots 
j_{k_{n}}^{(m)}(a_{n}),
$$
where $a_{l} \in {\cal G}_{k_{l}}$, $l=1, \ldots , n$.\\
\indent{\par}
{\sc Remark 1. }
We can also write the defining relations of Definition 2.1
in the following way:
$$
j_{1}^{(m)}(a)=\sum_{k=1}^{m} 
j_{1,k}^{(m)}(a),\;\;\;
j_{2}^{(m)}(a)=\sum_{k=1}^{m} 
j_{2,k}^{(m)}(a)
$$
where
$$
j_{1,k}^{(m)}(a)=i_{k,m}(a)\otimes (t_{[k,m]}-t_{[k-1,m]}),
$$
$$
j_{2,k}^{(m)}(b)=(t_{[k,m]}-t_{[k-1,m]})\otimes i_{k,m}(b),
$$
where $a \in {\cal G}_{1}$, $b\in {\cal G}_{2}$ and
we understand that $t_{[0,m]}=0$. It turns out that both
ways of writing Definition 2.1 (and its generalizations
introduced later) are useful, the first one -- for the
*-bialgebra construction, the second one -- for recurrence
relations. We will use them interchangably.\\
\indent{\par}
{\sc Remark 2.} The following notation will also be used:
$$
{\cal A}^{(m)}=j^{(m)}({\cal A}_{1}*{\cal A}_{2})\;\;\;
{\rm and}\;\;\;
\Phi^{(m)}=\widetilde{\Phi}^{(m)}|_{{\cal A}^{(m)}}.
$$
Moreover, the state $\Phi^{(m)}\circ j^{(m)}$ on 
${\cal A}_{1}*{\cal A}_{2}$ will be called the $m$-{\it free
product state}.\\
\indent{\par}
{\sc Remark 3.}
By $L^{(m)}$ we denote the two-sided ideal in $\widetilde{\cal
A}^{(m)}$ generated by $i_{k,2m}(t({\bf 1}-t))$, $1\leq k\leq 2m$.
It is easy to see that $L^{(m)}\subset
{\rm ker}\;\widetilde{\Phi}^{(m)}$. \\
\indent{\par}
{\sc Proposition 2.2.}
{\it Let $w,v$ are non-empty words in ${\cal A}_{1}, {\cal
A}_{2}$, respectively. Then}
$$
j_{1}^{(m)}(w)=\sum_{k=1}^{m}j_{1,k}^{(m)}(w)\;\; ({\rm
mod}\;\;L^{(m)}),\;\;\;
j_{2}^{(m)}(v)=\sum_{k=1}^{m}j_{2,k}^{(m)}(w)\;\; ({\rm
mod}\;\;L^{(m)}),
$$
{\it where}
$$
j_{1,k}^{(m)}(w)=i_{k,m}(w)\otimes (t_{[k,m]}-t_{[k-1,m]}),
$$
$$
j_{2,k}^{(m)}(v)=(t_{[k,m]}-t_{[k-1,m]})\otimes i_{k,m}(v).
$$
{\it Proof.} 
Let $a, a' \in {\cal G}_{1}$. We have
$$
j_{1}^{(m)}(a)j_{1}^{(m)}(a')=
\sum_{k,l=1}^{m}j_{1,k}^{(m)}(a)j_{1,l}^{(m)}(a').
$$
If $1<k<l$, then we obtain
$$
j_{1,k}^{(m)}(a)j_{1,l}^{(m)}(a')=
i_{k,m}(a)i_{l,m}(a')\otimes
(t_{[k,m]}-t_{[k-1,m]})(t_{[l,m]}-t_{[l-1,m]}) 
$$
$$
=i_{k,m}(a)i_{l,m}(a')\otimes
(I_{k-2}\otimes ({\bf 1}-t) \otimes t^{\otimes (m-k+1)})
(I_{l-2}\otimes ({\bf 1}-t) \otimes t^{\otimes (m-l+1)})
$$
$$
=i_{k,m}(a)i_{l,m}(a')\otimes
I_{k-2}\otimes ({\bf 1}-t) \otimes t^{\otimes (l-k-1)}\otimes
t({\bf 1}-t)
\otimes (t^{2})^{\otimes (m-l+1)} =0 \;\;({\rm mod}\;L^{(m)}).
$$
If $1=k < l$, then a similar analysis leads to 
$$
j_{1,1}^{(m)}(a)j_{1,l}^{(m)}(a')=
i_{1,m}(a)i_{l,m}(a')\otimes
t_{[1,m]}
(t_{[l,m]}-t_{[l-1,m]})
$$
$$
=i_{1,m}(a)i_{l,m}(a')\otimes
t^{\otimes m}(I_{l-2}\otimes ({\bf 1}-t) \otimes t^{\otimes
(m-l+1)})
$$
$$
= i_{1,m}(a)i_{l,m}(a')\otimes
t^{\otimes l-2}\otimes t({\bf 1}-t) \otimes (t^{2})^{\otimes
(m-l+1)}=0 
\;\;({\rm mod}\;L^{(m)}).
$$
Due to commutations, the case $k>l$ does not have to be
considered. Now,
$$
j_{1,k}^{(m)}(a)j_{1,k}^{(m)}(a')=i_{k,m}(aa')\otimes 
(t_{[k,m]}-t_{[k-1,m]})(t_{[k,m]}-t_{[k-1,m]}) 
$$
$$
=i_{k,m}(aa')\otimes 
(I_{k-2}\otimes ({\bf 1}-t)^{2} \otimes (t^{2})^{\otimes
(m-k+1)})=j_{1,k}^{(m)}(aa') \;\; ({\rm mod}\;L^{m}).
$$
This reasoning can now be extended to a product of $n$ generators.
The proof for $j_{2}^{(m)}(v)$ is similar.
\hfill $\Box$\\
\indent{\par}
{\sc Remark.}
Note that if we considered not free *-algebras, but
*-algebras, for which in the decomposition ${\cal A}_{k}={\cal A}_{k}^{0}\oplus
{\bf C}{\bf 1}$, ${\cal A}_{k}^{0}$ is a *-subalgebra of ${\cal A}_{k}$,
then we could obtain such algebras from the associated free *-algebras
by considering relations that do not involve the units. But from the above
proposition 
it is easy to see that such relations are preserved by $j^{(m)}_{k}$,
$k=1,2$ (modulo $L^{(m)}$). 
Therefore, our construction will also be valid for such unital
*-algebras.

Before we consider the general case, we look at 
the simplest case first, i.e. $m=1$. Then 
$$
\widetilde{{\cal A}}^{(1)}=\widetilde{{\cal A}}_{1}\otimes \widetilde{{\cal
A}}_{2}, \;\;\;
\widetilde{\Phi}^{(1)}=\widetilde{\phi}_{1}\otimes \widetilde{\phi}_{2}
$$
and
$$
j^{(1)}_{1}(a)=a\otimes t,\;\; j^{(1)}_{2}(b)=t\otimes b,
$$
where $a, b$ are generators of ${\cal A}_{1}, {\cal A}_{2}$,
respectively. Thus, if $w,v$ are nonempty words in ${\cal A}_{1}$, ${\cal
A}_{2}$, respectively, then
$$
j^{(1)}_{1}(w)=w\otimes t^{l(w)}=
w\otimes t\;\;({\rm mod}\ L^{(1)})
$$
$$
j^{(1)}_{2}(v)=
t^{l(v)}\otimes v\;\;=t\otimes v\;\;
({\rm mod}\;L^{(1)})
$$
where $l(w), l(v)$ are the lenghts of words $w,v$, respectively.

We obtain for $m=1$ the Boolean factorization law:
$$
\widetilde{\Phi}^{(1)}
\left(
j_{k_{1}}^{(1)}(w_{1})\ldots j_{k_{n}}^{(1)}(w_{n})
\right)
=\phi_{k_{1}}(w_{1})\ldots \phi_{k_{n}}(w_{n})
$$
with $\widetilde{\Phi}^{(1)}({\bf 1})=1$. Thus we can write
$$
\widetilde{\Phi}^{(1)}\circ j^{(1)}\equiv \phi_{1}*_{B}\phi_{2},
$$
where $\phi_{1}*_{B}\phi_{2}$ denotes the Boolean product of
$\phi_{1}$ and $\phi_{2}$. Thus the Boolean model is associated wih
$1$-freeness. This terminology can be justified
by means of the following argument: $\widetilde{\Phi}^{(1)}\circ j^{(1)}$ 
agrees with the conditionally free product on words 
$w_{1}w_{2}$.
This is the trivial case but simple enough to see how the ideas
of our approach developed. In the sequel we will construct
succesive approximations of (conditional) freeness, using tensor
products of higher orders. 

Before we go on, let us write down a simple result for the
Boolean case which will be used later.\\
\indent{\par}
{\sc Proposition 2.3.}
{\it The following factorization property holds:}
$$
\widetilde{\Phi}^{(1)}\left(
(j_{k_{1}}^{(1)}(w_{1})-d_{k_{1}}^{(1)}(w_{1}))
 \ldots
(j_{k_{n}}^{(1)}(w_{n})-d_{k_{n}}^{(1)}(w_{n}))
\right)
$$
$$
=
\left( 
\phi_{k_{1}}(w_{1})-\psi_{k_{1}}(w_{1})
\right)
\ldots
\left(
\phi_{k_{n}}(w_{n})-\psi_{k_{n}}(w_{n})
\right)
$$
{\it where $w_{1}, \ldots, w_{n}$ are non-empty words from
${\cal A}_{k_{1}}, \ldots , {\cal A}_{k_{n}}$, 
$k_{1}\neq k_{2}\neq \ldots \neq k_{n}$, and }
$$
d_{1}^{(1)}(w)=\psi_{1}(w) \otimes t,\;\;\;
d_{2}^{(1)}(w)=t\otimes \psi_{2}(w).
$$
{\it Proof.}
This property follows directly from the fact that $t$ plays in 
both $\widetilde{{\cal A}}_{1}$ and $\widetilde{{\cal A}}_{2}$ the
role of a separator of non-empty words from ${\cal
A}_{1}$ and ${\cal A}_{2}$, respectively. 
\hfill $\Box$

The next order of freeness will be associated with the double
tensor product
$$
\widetilde{{\cal A}}^{(2)}=\widetilde{{\cal A}}_{1}\otimes \widetilde{{\cal
A}}_{1}\otimes \widetilde{{\cal A}}_{2} \otimes \widetilde{{\cal A}}_{2}
$$
and the double tensor product state
$$
\widetilde{\Phi}^{(2)}=\widetilde{\phi}_{1}\otimes
\widetilde{\psi}_{1}\otimes \widetilde{\phi}_{2} \otimes
\widetilde{\psi}_{2} 
$$
with the *-homomorphism $j^{(2)}$ defined by
$$
j_{1}^{(2)}(a)=i_{1}(a)\otimes t_{[1,2]} 
+i_{2}(a)\otimes (t_{[2,2]}-t_{[1,2]})
$$
$$
\equiv
a\otimes {\bf 1}\otimes t \otimes t + {\bf 1}\otimes a
\otimes ({\bf 1}-t) \otimes t
$$
$$
j_{2}^{(2)}(b)
=t_{[1,2]}\otimes i_{1}(b)+(t_{[2,2]}-t_{[1,2]})\otimes i_{2}(b)
$$
$$
\equiv t \otimes t\otimes b \otimes {\bf 1} +
({\bf 1}-t) \otimes t \otimes {\bf 1} \otimes b
$$
for generators $a\in {\cal G}_{1},\;b\in {\cal G}_{2}$,
respectively. Let us present two examples.\\
\indent{\par}
{\sc Example 1.}
Let $a_{1},a_{2}\in {\cal G}_{1}$, $b\in {\cal G}_{2}$.
Then
$$
j_{1}^{(2)}(a_{1})j_{2}^{(2)}(b)j_{1}^{(2)}(a_{2})=
\left(
a_{1}\otimes {\bf 1}\otimes t \otimes t +
{\bf 1}\otimes a_{1}\otimes ({\bf 1}-t)\otimes t
\right)
$$
$$
\times
\left(
t\otimes t\otimes b \otimes {\bf 1} +({\bf 1}-t)\otimes t\otimes
{\bf 1} \otimes b
\right)
\left(
a_{2}\otimes {\bf 1}\otimes t \otimes t +
{\bf 1}\otimes a_{2}\otimes ({\bf 1}-t)\otimes t
\right)
$$
$$
=\left(
a_{1}\otimes {\bf 1}\otimes t \otimes t
\right)
\left(
t\otimes t\otimes b \otimes {\bf 1} +
({\bf 1}-t)\otimes t\otimes {\bf 1} \otimes b
\right)
\left(
a_{2}\otimes {\bf 1}\otimes t \otimes t 
\right)
\;\;({\rm mod}\;L^{(2)}).
$$
Therefore,
$$
\Phi^{(2)}
\left(
j_{1}^{(2)}(a_{1})j_{2}^{(2)}(b)j_{1}^{(2)}(a_{2})
\right)
=
\widetilde{\Phi}^{(2)}(a_{1}ta_{2}\otimes t \otimes tbt \otimes t)
$$
$$
+
\widetilde{\Phi}^{(2)}(a_{1}a_{2}\otimes t \otimes t\otimes tbt)-
\widetilde{\Phi}^{(2)}(a_{1}ta_{2}\otimes t\otimes t\otimes tbt)
$$
$$
=\phi_{2}(a_{1})\phi_{2}(a_{2})\phi_{1}(b)+
\phi_{2}(a_{1}a_{2})\psi_{1}(b) -
\phi_{2}(a_{1})\phi_{2}(a_{2})\psi_{1}(b).
$$
\indent{\par}
{\sc Example 2.}
Let $a_{1},a_{2}\in {\cal G}_{1}$, $b_{1}, b_{2}\in {\cal
G}_{2}$. Then
$$
j_{1}^{(2)}(a_{1})j_{2}^{(2)}(b_{1})j_{1}^{(2)}(a_{2})j_{2}^{(2)}(b_{2})=
\left(
a_{1}t\otimes t \otimes tb_{1}\otimes t+
a_{1}({\bf 1}-t)\otimes t \otimes t \otimes tb_{1}
\right)
$$
$$
\times
\left(
a_{2}t\otimes t \otimes tb_{2}\otimes t+
t\otimes a_{2}t\otimes ({\bf 1}-t)b_{2}\otimes t
\right) \;\;\; ({\rm mod}\;L^{(2)})
$$
$$
= a_{1}ta_{2}t\otimes t\otimes  tb_{1}tb_{2}\otimes t
+
a_{1}t\otimes ta_{2}t\otimes
\left[
 tb_{1}b_{2}\otimes t - tb_{1}tb_{2} \otimes t
\right]
$$
$$
+
\left[
 a_{1}a_{2}t \otimes t -  a_{1}ta_{2}t \otimes t
\right]
\otimes
 tb_{2} \otimes tb_{1}t \;\;\;({\rm mod} \;L^{(2)}).
$$
Therefore, we obtain
$$
{\Phi}^{(2)}\left(
j_{1}^{(2)}(a_{1})j_{2}^{(2)}(b_{1})j_{1}^{(2)}(a_{2})j_{2}^{(2)}(b_{2})
\right)
=
\phi_{1}(a_{1})\phi_{1}(a_{2})
\phi_{2}(b_{1})\phi_{2}(b_{2})
$$
$$
+\phi_{1}(a_{1})\psi_{1}(a_{2})
\left[
\phi_{2}(b_{1}b_{2})-
\phi_{2}(b_{1})\phi_{2}(b_{2})
\right]
+
\left[
\phi_{1}(a_{1}a_{2})-\phi_{1}(a_{1})\phi_{1}(a_{2})
\right] 
\phi_{2}(b_{2})\psi_{2}(b_{1}).
$$
In both examples we obtain the same expressions as if we
calculated $*_{l\in \{1,2\}}(\phi_{l},\psi_{l})$
acting on $a_{1}ba_{2}$ and $a_{1}b_{1}a_{2}b_{2}$, respectively. 
\\[5pt]
\newpage
\begin{center}
{\sc 3. Factorization Lemmas}
\end{center}
In this section we will derive some factorization lemmas that
will be needed in the proofs of the main theorems in Section 4. 

Let us define the following ``condition'' maps 
${\Psi}^{(m)}$:
$$
{\Psi}^{(m)}={\rm id}^{\otimes (m-1)}\otimes
\widetilde{\psi}_{1} \otimes {\rm id}^{\otimes (m-1)}\otimes
\widetilde{\psi}_{2}.
$$
Note that 
$$
\widetilde{\Phi}^{(m)}=\widetilde{\Phi}^{(m-1)}\circ {\Psi}^{(m)}
=\widetilde{\Phi}^{(m-2)}\circ
{\Psi}^{(m-1)} \circ {\Psi}^{(m)} = \ldots =
\widetilde{\Phi}^{(1)}\circ {\Psi}^{(2)}\circ \ldots 
\circ {\Psi}^{(m)}.
$$
\indent{\par}
{\sc Proposition 3.0.}
{\it We can write}
$$
\left(
{\Psi}^{(m)}\circ j_{1}^{(m)}
\right)
(w)
=j_{1}^{(m-1)}(w) +g_{1}^{(m-1)}(w) \;\;\; 
({\rm mod}\; L^{(m-1)})
$$
$$
\left(
{\Psi}^{(m)}\circ j_{2}^{(m)}
\right)
(v)
=j_{2}^{(m-1)}(v) +g_{2}^{(m-1)}(v) \;\;\; 
({\rm mod}\; L^{(m-1)})
$$
{\it where}
$$
g_{1}^{(m-1)}(w)=\psi_{1}(w)[I_{m-1}\otimes I_{m-2}\otimes 
({\bf 1}-t)],
$$
$$
g_{2}^{(m-1)}(v)=\psi_{2}(v)[I_{m-2}\otimes ({\bf 1}-t)\otimes I_{m-1}],
$$
{\it and $w,v$ are non-empty words in ${\cal A}_{1}, {\cal
A}_{2}$, respectively.}\\[5pt]
{\it Proof.} It is an immediate consequence of Proposition 2.2. 
\hfill $\Box$.\\
\indent{\par}
{\sc Proposition 3.1.}
{\it The unital *-homomorphisms $j^{(m)}$ preserve the marginal
laws, i.e.}
$$
\widetilde{\Phi}^{(m)} \circ j_{1}^{(m)} 
= \phi_{1},\;\;\;
\widetilde{\Phi}^{(m)} \circ j_{2}^{(m)}
= \phi_{2}.
$$
{\it Proof.} 
If $m=1$, it is obvious. For $m>1$, we obtain the result using
the induction argument.
Clearly, $\left(\Phi^{(m)} \circ
j_{k}^{(m)}\right)({\bf 1})=1=\phi_{k}({\bf 1})$, 
since $\phi_{l}, \psi_{l}$ are
states on ${\cal A}_{l}$, $l=1,2$. Thus, assume that
$w=a_{1}\ldots a_{n}$ is a non-empty word in ${\cal A}_{1}$.
Using Proposition 3.0, we obtain
$$
\left(
\widetilde{\Phi}^{(m)}\circ
j_{1}^{(m)}
\right)
(w)
=\widetilde{\Phi}^{(m-1)}\circ {\Psi}^{(m)} 
\left(
j_{1}^{(m-1)}(w)+g_{1}^{(m-1)}(w)
\right)
=\widetilde{\Phi}^{(m-1)}
\left(
j_{1}^{(m-1)}(w)
\right).
$$
\hfill $\Box$

In the expressions for $j_{k}^{(m)}(w)$, $k=1,2$, there is always
exactly one term with one separator $t$.
Since it will have to be subtracted from $j_{k}^{(m)}(w)$, we introduce a new
notation. Thus, for words $w,v$ in 
${\cal A}_{1}, {\cal A}_{2}$, respectively, let
$$
d_{1}^{(m)}(w)=i_{m,m}(w)\otimes t_{[m,m]}\equiv
I_{m-1}\otimes w\otimes I_{m-1} \otimes t,
$$
$$
d_{2}^{(m)}(v)=t_{[m,m]}\otimes i_{m,m}(v)\equiv
I_{m-1}\otimes t\otimes I_{m-1}\otimes v .
$$
Let us also note that
$$
({\Psi}^{(m)} \circ d_{1}^{(m)})(w)=
\psi_{1}(w)[I_{m-1}\otimes I_{m-1}],\;\;\;
({\Psi}^{(m)} \circ d_{2}^{(m)})(v)=
\psi_{2}(v)[I_{m-1}\otimes I_{m-1}].
$$
for words $w,v$ in ${\cal A}_{1}, {\cal A}_{2}$, respectively.
Moreover,
$$
{\Psi}^{(m)}(j_{k}^{(m)}(w)-d_{k}^{(m)}(w))
=j_{k}^{(m-1)}(w)-h_{k}^{(m-1)}(w) \;\;({\rm mod}\;L^{(m-1)})
$$
{\it where $k=1,2$ and }
$$h_{1}^{(m-1)}(w)=\psi_{1}(w)[I_{m-1}\otimes I_{m-2}\otimes t],\;\;\;
h_{2}^{(m-1)}(w)=\psi_{2}(w)[I_{m-2}\otimes t \otimes I_{m-1}]
$$
The main purpose of introducing the separator $t$ is to
obtain some factorizations of correlations. We present two easy
lemmas. \\
\indent{\par}
{\sc Lemma 3.2.}
{\it 
Let $w_{1}, \ldots, w_{n}$ be non-empty words in ${\cal A}_{k_{1}},
\ldots ,{\cal A}_{k_{n}}$, respectively, where $k_{1}\neq
k_{2}\neq \ldots \neq k_{n}$. Then ${\Psi}^{(m)}$ exhibits the following
multiplicative property:}
$$
\left(
{\Psi}^{(m)} \circ j^{(m)}
\right)
\left(
w_{1}\ldots w_{n}
\right)
= \left(
{\Psi}^{(m)} \circ j_{k_{1}}^{(m)}
\right)
\left(
w_{1}
\right)
\ldots
\left(
{\Psi}^{(m)} \circ j_{k_{n}}^{(m)}
\right)
\left(
w_{n}
\right)\;\;({\rm mod}\;L^{(m-1)})
$$
{\it Proof.}
The only thing to show is that all words from ${\cal A}_{1}$ 
appearing at the $m$-th site and
all words from ${\cal A}_{2}$ appearing at the $2m$-th site are
separated by $t$. But that immediately follows from Proposition
2.2 since each summand of $j_{1}^{(m)}(w)$, $w\in {\cal A}_{1}$,
has $t$ at the $2m$-th site and each summand of
$j_{2}^{(m)}(v)$, $v\in {\cal A}_{2}$, has $t$ at the $m$-th site.
\hfill
$\Box$\\
\indent{\par}
{\sc Lemma 3.3.}
{\it Let $w_{1}, \ldots, w_{p}$ be non-empty words in ${\cal A}_{k_{1}},
\ldots, {\cal A}_{k_{n}}$, respectively, where $k_{1}, \ldots ,
k_{n}\in \{1,2\}$ and $k_{1}\neq k_{2}\neq \ldots \neq k_{n}$.
Then} 
$$
\widetilde{\Phi}^{(m)} 
\left[
\left(
j_{k_{1}}^{(m)}(w_{1})-d_{k_{1}}^{(m)}(w_{1})
\right)
\ldots
\left(
j_{k_{n}}^{(m)}(w_{n})-d_{k_{n}}^{(m)}(w_{n})
\right)
\right]
$$
$$
=
\left( 
\phi_{k_{1}}(w_{1})-\psi_{k_{1}}(w_{1})
\right)
\ldots
\left(
\phi_{k_{n}}(w_{n})-\psi_{k_{n}}(w_{n})
\right).
$$
{\it Proof.}
The case $m=1$ is covered by Proposition 2.1. The general case
follows from the induction argument.
Note that $\widetilde{\Phi}^{(m)}=\widetilde{\Phi}^{(m-1)}\circ \Psi^{(m)}$ and use
Lemma 3.2 to obtain
$$
\widetilde{\Phi}^{(m)} 
\left[
\left(
j_{k_{1}}^{(m)}(w_{1})-d_{k_{1}}^{(m)}(w_{1})
\right)
\ldots
\left(
j_{k_{n}}^{(m)}(w_{n})-d_{k_{n}}^{(m)}(w_{n})
\right)
\right]
$$
$$
=\widetilde{\Phi}^{(m-1)} 
\left[
\left(
j_{k_{1}}^{(m-1)}(w_{1})-h_{k_{1}}^{(m-1)}(w_{1})
\right)
\ldots
\left(
j_{k_{n}}^{(m-1)}(w_{n})-h_{k_{n}}^{(m-1)}(w_{n})
\right)
\right]
$$
Write $\widetilde{\Phi}^{(m-1)}=\widetilde{\Phi}^{(m-2)}\circ
\Psi^{(m-1)}$ and, since each $h_{1}^{(m-1)}(w_{l})$ has $t$ at the
$2m-2$-th site and each $h_{2}^{(m-1)}(w_{k})$ has $t$ at the
$m-1$-th site,
$\Psi^{(m-1)}$ is multiplicative also on the products of the
above type, hence we obtain
$$
\widetilde{\Phi}^{(m-2)}
\left[
\Psi^{(m-1)}
\left(
j_{k_{1}}^{(m-1)}(w_{1})-h_{k_{1}}^{(m-1)}(w_{1})
\right)
\ldots
\Psi^{(m-1)}
\left(
j_{k_{n}}^{(m-1)}(w_{n})-h_{k_{1}}^{(m-1)}(w_{n})
\right)
\right]
$$
$$
=
\widetilde{\Phi}^{(m-2)}
\left[
\Psi^{(m-1)}
\left(
j_{k_{1}}^{(m-1)}(w_{1})-d_{k_{1}}^{(m-1)}(w_{1})
\right)
\ldots
\Psi^{(m-1)}
\left(
j_{k_{n}}^{(m-1)}(w_{n})-d_{k_{1}}^{(m-1)}(w_{n})
\right)
\right]
$$
$$
=
\widetilde{\Phi}^{(m-1)}
\left[
\left(
j_{k_{1}}^{(m-1)}(w_{1})-d_{k_{1}}^{(m-1)}(w_{1})
\right)
\ldots
\left(
j_{k_{n}}^{(m-1)}(w_{n})-d_{k_{1}}^{(m-1)}(w_{n})
\right)
\right]
$$
where (in the first equation) we used
$$
h_{1}^{(m-1)}(w)-
d_{1}^{(m-1)}(w)\in \;{\rm ker}{\Psi}^{(m-1)}
$$
$$
h_{2}^{(m-1)}(v)-
d_{2}^{(m-1)}(v)\in \;{\rm ker}{\Psi}^{(m-1)}
$$
for words $w,v$ in ${\cal A}_{1}$, ${\cal A}_{2}$, respectively.
This finishes the proof.
\hfill
$\Box$\\[5pt]
\begin{center}
{\sc 4. Main Theorems}
\end{center}
We are ready to state our main result which says that
the $m$-free product state agrees with the (conditionally) 
free product state on word products of not more than $2m$ words.
The proof of that fact will be carried out in two steps. First,
we show that $m$-freeness agrees with conditional freeness for
products of at most $m+1$ words (Theorem 4.0). Then, we will
improve that result in Theorem 4.1 and prove that in fact $m+1$
can be replaced by $2m$.
As a corollary we obtain pointwise
convergence of the $m$-free product 
states $\Phi^{(m)}\circ j^{(m)}$ to the
conditionally free product state. Thus, the conditionally free
case, in particular the free case, is obtained as a limit of
$m$-fold tensor product constructions.\\
\indent{\par}
{\sc Theorem 4.0.}
{\it Let
$\widetilde {\Phi}^{(m)}=
\widetilde{\Phi}_{1}^{(m)}\otimes \widetilde{\Phi}_{2}^{(m)}$, 
$\widetilde{\Phi}_{l}^{(m)}=
\widetilde{\phi}_{l}\otimes \widetilde{\psi}_{l}^{(m-1)}$, $l=1,2$, where
$\widetilde{\phi}_{l}$, $\widetilde{\psi}_{l}$ are Boolean extensions of
states $\phi_{l}$, $\psi_{l}$ on unital free *-algebras, 
${\cal A}_{i}$, $i=1,2$. Then, if $n\leq m+1$, then 
$\widetilde{\Phi}^{(m)} \circ j^{(m)}$
agrees with the conditionally free product
$*_{i\in \{1,2\}}(\phi_{i},\psi_{i})$
on word products $w_{1}\ldots w_{n}$, where 
$w_{1},\ldots ,w_{n}\in {\cal A}_{k_{1}}, \ldots ,
{\cal A}_{k_{n}}$, respectively, and 
$k_{1}\neq k_{2}\neq \ldots \neq k_{n}$}\\[5pt]
{\it Proof.}
If $m=1$, then the result is trivial (Boolean case).
So let us proceed with the induction. 
We have $(\Phi^{(m)} \circ j^{(m)})({\bf 1})=\Phi^{(m)}(I_{m}\otimes
I_{m})=1$. Now,
$$
{\Phi}^{(m)} 
\left(
j^{(m)}_{k_{1}}(w_{1})\ldots 
j^{(m)}_{k_{n}}(w_{n})
\right)
$$
$$
=
\widetilde{\Phi}^{(m)}
\left(
\left(
j^{(m)}_{k_{1}}(w_{1})-d^{(m)}_{k_{1}}(w_{1})
\right)
\ldots 
\left(
j^{(m)}_{k_{n}}(w_{n})-d^{(m)}_{k_{n}}(w_{n})
\right)
\right)
$$
$$
+
\sum_{i}
\widetilde{\Phi}^{(m)} 
\left(
j^{(m)}_{k_{1}}(w_{1})\ldots
d^{(m)}_{k_{i}}(w_{i})\ldots 
j^{(m)}_{k_{n}}(w_{n})
\right)
$$
$$
-
\sum_{i<l}
\widetilde{\Phi}^{(m)} 
\left(
j^{(m)}_{k_{1}}(w_{1})\ldots 
d^{(m)}_{k_{i}}(w_{i})\ldots
d^{(m)}_{k_{l}}(w_{l})\ldots
j^{(m)}_{k_{n}}(w_{n})
\right)
$$
$$
+ \ldots - (-1)^{n}
\widetilde{\Phi}^{(m)} 
\left(
d^{(m)}_{k_{1}}(w_{1})\ldots 
d^{(m)}_{k_{n}}(w_{n})
\right).
$$
Note that the above recurrence relation looks like the corresponding
one for the conditionally free case, except that instead of numbers
we have $d^{(m)}_{k_{i}}(w)$'s 
replacing $j^{(m)}_{k_{i}}(w)$'s at one, two, or more places.
We invoke Lemma 3.2 to conclude that the first term on the right
-hand side is identical to the conditionally free case. Thus
what we need to prove is that in the remaining ones 
$d^{(m)}_{k_{i}}(w)$'s
indeed behave like numbers.
Thus the proof reduces to proving the following claim.\\[5pt]
{\sc Claim:}
$$
\widetilde{\Phi}^{(m)} 
\left(
j^{(m)}_{k_{1}}(w_{1})\ldots
d^{(m)}_{k_{i(1)}}(w_{i(1)})\ldots 
d^{(m)}_{k_{i(l)}}(w_{i(l)})\ldots 
j^{(m)}_{k_{n}}(w_{n})
\right)
$$
$$
=*_{i\in \{1,2\}}(\phi_{i}, \psi_{i})
(w_{1}\ldots \breve{w}_{i(1)}\ldots  \breve{w}_{i(l)} \ldots w_{n})
\psi_{k_{i(1)}}(w_{i(1)})\ldots \psi_{k_{i(l)}}(w_{i(l)}),
$$
for $n\leq m+1$, where by $\breve{}$ 
we understand that the words with indices
$i(1), \ldots , i(l)$, $1\leq l\leq n$, are omitted. 
Note that the claim says that 
the operators $d^{(m)}_{k_{i}}(w)$ do behave like constants 
locally, i.e. if the correlation is not too long (for now, $n\leq
m+1$) and that is the reason why we only have local freeness for
fixed $m$. 

The claim will be proved by induction. It trivially holds for
$m=1$ and $n\leq 2$, so assume that it holds for $m-1$.
In particular, this inductive assumption implies that 
$\widetilde{\Phi}^{(m-1)}\circ j^{(m-1)}$ agrees with 
$*_{i\in \{1,2\}}(\phi_{i},\psi_{i})$ on products of 
$n\leq m$ words (Lemma 3.3 is used and the above recurrence
relation for $m-1$).

Using $\widetilde{\Phi}^{(m)}=\widetilde{\Phi}^{(m-1)}\circ \Psi^{(m)}$, we obtain
$$
\widetilde{\Phi}^{(m)} 
\left(
j^{(m)}_{k_{1}}(w_{1})\ldots
d^{(m)}_{k_{i(1)}}(w_{i(1)})\ldots 
d^{(m)}_{k_{i(l)}}(w_{i(l)})\ldots 
j^{(m)}_{k_{n}}(w_{n}) 
\right)
$$
$$
=
\widetilde{\Phi}^{(m-1)} 
\left(
\Psi^{(m)}
(
j^{(m)}_{k_{1}}(w_{1})
)
\ldots
\Psi^{(m)}
(
d^{(m)}_{k_{i(1)}}(w_{i(1)})
)
\ldots 
\Psi^{(m)}
(
d^{(m)}_{k_{i(l)}}(w_{i(l)})
)
\ldots 
\Psi^{(m)}
(
j^{(m)}_{k_{n}}(w_{n})
)
\right)
$$
$$
=
\psi_{k_{i(1)}}(w_{i(1)})\ldots \psi_{k_{i(l)}}(w_{i(l)})
\widetilde{\Phi}^{(m-1)} 
\left(
\prod_{i\in [n]\setminus \{i(1), \ldots , i(l)\}}
(
j^{(m-1)}_{k_{i}}(w_{i}) +g^{(m-1)}_{k_{i}}(w_{i})
)
\right).
$$
We used the muliplicativity of $\Psi^{(m)}$ (cf. Lemma 3.2),
which still holds when some of the
$j^{(m)}_{k_{i}}(w_{i})$'s are replaced by $d^{(m)}_{k_{i}}(w_{i})$'s
($d^{(m)}_{1}(w)$ and $d^{(m)}_{2}(v)$ have $t$ at
the $2m$-th and $m$-th tensor sites, respectively, and $k_{1}\neq
k_{2}\neq \ldots \neq k_{n}$). 
We then used Proposition 3.0 to get the second equation.

It is enough to show that
$$
\widetilde{\Phi}^{(m-1)} 
\left(
j^{(m-1)}_{k_{1}}(w_{1})\ldots
g^{(m-1)}_{k_{p(1)}}(w_{p(1)})\ldots 
g^{(m-1)}_{k_{p(r)}}(w_{p(r)})\ldots 
j^{(m-1)}_{k_{n}}(w_{n}) 
\right)=0
$$
for $n\leq m$ and arbitrary $k_{1}, \ldots, k_{n}$
(note that we pulled out at least one
$d_{k_{i}}^{(m)}(w_{i})$ above, so the number of factors got smaller). 

Let us now make some observations which will
reduce the number of cases that need to be considered. 
We refer to 
$g^{(m-1)}_{k_{p(1)}}(w_{p(1)}),\ldots ,
g^{(m-1)}_{k_{p(r)}}(w_{p(r)}) 
$
in the above formula, 
although, for simplicity, 
the indices $k_{p(1)}, \ldots , k_{p(r)}$ and $p(1), \ldots,
p(r)$ will not be used explicitly. Instead, we will refer to
generic $w,w'$ or $v$, nonempty words in ${\cal A}_{1}$, ${\cal A}_{2}$,
respectively.
Firstly,
note that each element of type $g_{1}^{(m-1)}(w)$ commutes with
$g_{2}^{(m-1)}(v)$, so if such elements stand next to
each other, we can regroup them in any way we want.
Secondly, without loss of generality we can replace 
$g_{k}^{(m-1)}(w)g_{k}^{(m-1)}(w')$ by $g_{k}^{(m-1)}(ww')$,
$k=1,2$ (this only changes the above expression by a constant).
Thirdly, it is enough to consider such configurations
in which each $g^{(m-1)}_{1}(w)$ is sourrounded by 
$j^{(m-1)}_{2}(w')$'s and $g^{(m-1)}_{2}(w)$ is sourrounded by
$j^{(m-1)}_{1}(w')$'s. For, if for instance we had
$g_{1}^{(m-1)}(w)j^{(m-1)}_{1}(w')$ or  
$j^{(m-1)}_{1}(w')g_{1}^{(m-1)}(w)$, then at the $2m-2$-th
tensor site we would get $t({\bf 1}-t)$, which is in $L^{(m-1)}$.
A similar argument shows that $g^{(m-1)}_{2}(w)$ must be
sourronded by $j^{(m-1)}_{1}(w')$'s. All this reduces the 
proof to configurations in which elements
$g_{k_{p(i)}}^{(m-1)}(w_{p(i)})$ are sourrounded by 
$j_{k_{p(i)-1}}^{(m-1)}(w_{p(i)-1})$ and 
$j_{k_{p(i)+1}}^{(m-1)}(w_{p(i)+1})$ with $k_{p(i)-1}\neq
k_{p(i)}\neq k_{p(i)+1}$. Thus it remains to tackle the
configurations of this type and show that their contribution
vanishes. 

For that purpose we will replace each $g_{1}^{(m-1)}(w)$ and
$g_{2}^{(m-1)}(v)$ by 
$$
\psi_{1}(w)I_{2m-2} -h_{1}^{(m-1)}(w),\;\;\;
\psi_{2}(v)I_{2m-2} -h_{2}^{(m-1)}(v),
$$
respectively.
Thus, to finally prove the claim it suffices to show that for $n\leq m$ we
have 
$$
\widetilde{\Phi}^{(m-1)} 
\left(
j^{(m-1)}_{k_{1}}(w_{1})\ldots
h^{(m-1)}_{k_{s(1)}}(w_{p(1)})\ldots 
h^{(m-1)}_{k_{s(u)}}(w_{p(r)})\ldots 
j^{(m-1)}_{k_{n}}(w_{n}) 
\right)
$$
$$
=
\psi_{k_{s(1)}}(w_{s(1)})\ldots \psi_{k_{s(u)}}(w_{p(r)})
\widetilde{\Phi}^{(m-1)} 
\left(
\prod_{i\in [n]\setminus \{s(1), \ldots , s(u)\}}
(
j^{(m-1)}_{k_{i}}(w_{i})
)
\right)
$$
for $k_{1}\neq k_{2}\neq \ldots \neq k_{n}$. 
In fact, one can take the consecutive indices different since
$j^{(m-1)}_{1}$, $j^{(m-1)}_{2}$ are multiplicative and the
configurations to which we reduced our proof had the property
$k_{1}\neq k_{2}\neq \ldots \neq k_{n}$.

We write again 
$\widetilde{\Phi}^{(m-1)}=\widetilde{\Phi}^{(m-2)}\circ \Psi^{(m-1)}$ 
and use the factorization property of $\Psi^{(m-1)}$ on the
left-hand side of the above equation since $h_{1}^{(m-1)}$ and
$h_{2}^{(m-1)}$ have a non-zero power of $t$ at the $2m-2$-th
and $m-1$-th tensor site, respectively. We obtain
$$
\widetilde{\Phi}^{(m-1)} 
\left(
j^{(m-1)}_{k_{1}}(w_{1})\ldots
h^{(m-1)}_{k_{s(1)}}(w_{s(1)})\ldots 
h^{(m-1)}_{k_{s(u)}}(w_{s(u)})\ldots 
j^{(m-1)}_{k_{n}}(w_{n}) 
\right)
$$
$$
=\widetilde{\Phi}^{(m-2)}
\left(
\Psi^{(m-1)}
(
j^{(m-1)}_{k_{1}}(w_{1})
)
\ldots
\Psi^{(m-1)}
(
h^{(m-1)}_{k_{s(1)}}(w_{s(1)})
)
\right.
$$
$$
\ldots 
\left.
\Psi^{(m-1)}
(
h^{(m-1)}_{k_{s(u)}}(w_{s(u)})
)
\ldots
\Psi^{(m-1)}
( 
j^{(m-1)}_{k_{n}}(w_{n}) 
)
\right)
$$
$$
=\widetilde{\Phi}^{(m-2)}
\left(
\Psi^{(m-1)}
(
j^{(m-1)}_{k_{1}}(w_{1})
)
\ldots
\Psi^{(m-1)}
(
d^{(m-1)}_{k_{s(1)}}(w_{s(1)})
)
\right.
$$
$$
\left.
\ldots 
\Psi^{(m-1)}
(
d^{(m-1)}_{k_{s(u)}}(w_{s(u)})
)
\ldots
\Psi^{(m-1)}
( 
j^{(m-1)}_{k_{n}}(w_{n}) 
)
\right)
$$
$$
=
\widetilde{\Phi}^{(m-1)} 
\left(
j^{(m-1)}_{k_{1}}(w_{1})\ldots
d^{(m-1)}_{k_{s(1)}}(w_{s(1)})\ldots 
d^{(m-1)}_{k_{s(u)}}(w_{s(u)})\ldots 
d^{(m-1)}_{k_{n}}(w_{n}) 
\right)
$$
$$
=
\psi_{k_{s(1)}}(w_{s(1)})\ldots \psi_{k_{s(u)}}(w_{s(u)})
\widetilde{\Phi}^{(m-1)} 
\left(
\prod_{i\in [n]\setminus \{s(1), \ldots , s(u)\}}
(
j^{(m-1)}_{k_{i}}(w_{i})
)
\right)
$$
where, in the last equation, we used the inductive assumption of
the claim. Thus, we have proved our claim. This also finishes 
the proof of the theorem.
\hfill $\Box$

It turns out that the result of the theorem can be improved.
Namely, one can show that $\widetilde{\Phi}^{(m)}\circ j^{(m)}$ agrees
with the conditionally free product of states on word products
$w_{1}\ldots w_{n}$ for $n\leq 2m$. This is not surprising since
already Example 2 in Section 2 showed that $\widetilde{\Phi}^{(2)}\circ
j^{(2)}$ suffices to calculate a four-point correlation.\\
\indent{\par}
{\sc Theorem 4.1.}
{\it $\widetilde{\Phi}^{(m)}\circ j^{(m)}$ agrees with the conditionally
free product on word products $w_{1}\ldots w_{n}$ for $n\leq 2m$, where
$w_{1}, \ldots , w_{n}$ 
$\in {\cal A}_{k_{1}}, \ldots , {\cal A}_{k_{n}}$,
respectively, and $k_{1}\neq k_{2}\neq \ldots  \neq k_{n}$.}
\\[5pt]
{\it Proof.}
Let $n\leq 2m$. We know from Theorem 4.0 that
$\widetilde{\Phi}^{(s)}\circ j^{(s)}$, $s=2m-1$, agrees with the conditionally
free product on word products $w_{1}\ldots w_{n}$. 
We will show that $\widetilde{\Phi}^{(s)}\circ j^{(s)}$ agrees in fact with
$\widetilde{\Phi}^{(m)}\circ j^{(m)}$ on such word products. 
By Proposition 2.2 we have
$$
\widetilde{\Phi}^{(s)}
\left(
j_{k_{1}}^{(s)}(w_{1})\ldots j_{k_{n}}^{(s)}(w_{n})
\right)
=\sum_{m_{1}=1}^{s}\ldots \sum_{m_{n}=1}^{s}
\widetilde{\Phi}^{(s)}
\left(
j_{k_{1},m_{1}}^{(s)}(w_{1})
\ldots 
j_{k_{n}, m_{n}}^{(s)}(w_{n})
\right).
$$
{\sc Claim.}
One obtains the following ``pyramidal'' formula:
$$
\widetilde{\Phi}^{(s)}
\left(
j_{k_{1}}^{(s)}(w_{1})\ldots j_{k_{n}}^{(s)}(w_{n})
\right)
$$
$$
=\sum_{m_{1}=1}^{1}
\sum_{m_{2}=1}^{2}
\ldots
\sum_{m_{n-1}=1}^{2}
\sum_{m_{n}=1}^{1}
\widetilde{\Phi}^{(s)}
\left(
j_{k_{1},m_{1}}^{(s)}(w_{1})
\ldots 
j_{k_{n}, m_{n}}^{(s)}(w_{n})
\right).
$$
Note that if $n=2k$ is even we obtain a ``pyramid'' of height
$k\leq m$ with a flat top, 
and if $n=2k-1$ is odd, then we get a ``pyramid'' of
height $k\leq m$ with a sharp top.

To prove our claim we assume that $k_{1}=1$. The proof for
$k_{1}=2$ is similar.
First, note that the only term that survives from the first summation
corresponds to $m_{1}=1$. The reason is simple. That is the only term
that does not have a compensator since it takes the form
$$
j_{1,1}^{(s)}(w_{1})=i_{1,s}(w_{1})\otimes t_{[1,s]} .
$$
The other ones look like
$$
j_{1,r}^{(s)}(w_{1})=i_{r,s}(w_{1})\otimes (t_{[r,s]}-t_{[r-1,s]})
$$
for $r>1$ and thus have $t-{\bf 1}$ at site $(2,r)$ which is not preceded
by any non-empty words of ${\cal A}_{2}$ and thus give zero.
By mirror reflection we can conclude that the same must happen
at the other end of the correlation. Thus, we obtain
$$
\widetilde{\Phi}^{(s)}
\left(
j_{k_{1}}^{(s)}(w_{1})\ldots j_{k_{n}}^{(s)}(w_{n})
\right)
$$
$$
=\sum_{m_{2}=1}^{s}\ldots 
\sum_{m_{n-1}=1}^{s}
\widetilde{\Phi}^{(s)}
\left(
j_{k_{1},1}^{(s)}(w_{1})
j_{k_{2},m_{2}}^{(s)}(w_{2})
\ldots 
j_{k_{n-1},m_{n-1}}^{(s)}(w_{n-1})
j_{k_{n}, 1}^{(s)}(w_{n})
\right).
$$
Suppose that we have already reduced our expresion to the
following form
$$
\widetilde{\Phi}^{(s)}
\left(
j_{k_{1}}^{(s)}(w_{1})\ldots j_{k_{n}}^{(s)}(w_{n})
\right)
$$
$$
=\sum_{m_{1}=1}^{1}
\ldots
\sum_{m_{l}=1}^{l}
\sum_{m_{l+1}=1}^{s}
\ldots 
\sum_{m_{n-l}=1}^{s}
\sum_{m_{n-l+1}=1}^{l}
\ldots
\sum_{m_{n}=1}^{1}
\widetilde{\Phi}^{(s)}
\left(
j_{k_{1},m_{1}}^{(s)}(w_{1})
\ldots 
j_{k_{n},m_{n}}^{(s)}(w_{n})
\right).
$$
To fix attention, assume that $k_{l+1}=1$. 
We will show that the terms in which $j_{k_{l+1},
r}^{(s)}(w_{l+1})$ appears for $r>l+1$, give vanishing contribution.
Such a term produces
$$
i_{1,1}^{(s)}(w_{1})(t_{[m_{2},s]}-t_{[m_{2}-1,s]})\ldots 
(t_{[m_{l},s]}-t_{[m_{l}-1,s]})i_{r,s}(w_{l+1})\ldots
$$
in the place reserved for $\widetilde{{\cal A}}_{1}^{\otimes s}$
and
$$
t_{[1,s]}i_{m_{2},s}(w_{2})\ldots
i_{m_{l},s}(w_{l})(t_{[r,s]}-t_{[r-1,s]})\ldots
$$
in the place reserved for $\widetilde{{\cal A}}_{2}^{\otimes s}$.
The second expression is crucial. Namely, if $r>l+1$, then
the $t$'s produced by $j_{k_{l+1},r}^{(s)}(w_{l+1})$ appear at
sites greater than $l+1$. At these sites there are no words of
${\cal A}_{2}$ preceding the $t$'s. Therefore, the term with
$t_{[r,s]}$ is compensated by the term with $t_{[r-1,s]}$.
Again, the miror reflection gives a symmetric behavior on the
other side. The proof for $k_{l+1}=2$ is similar. This finishes
the proof of the claim.

Thus, we finally have to show that in order to perform
calculations for a pyramid of height $m$ one can replace $s$
by $m$, i.e.
$$
\sum_{m_{1}=1}^{1}
\sum_{m_{2}=1}^{2}
\ldots
\sum_{m_{n-1}=1}^{2}
\sum_{m_{n}=1}^{1}
\widetilde{\Phi}^{(s)}
\left(
j_{k_{1},m_{1}}^{(s)}(w_{1})
\ldots 
j_{k_{n}, m_{n}}^{(s)}(w_{n})
\right)
$$
$$
=\sum_{m_{1}=1}^{1}
\sum_{m_{2}=1}^{2}
\ldots
\sum_{m_{n-1}=1}^{2}
\sum_{m_{n}=1}^{1}
\widetilde{\Phi}^{(m)}
\left(
j_{k_{1},m_{1}}^{(m)}(w_{1})
\ldots 
j_{k_{n}, m_{n}}^{(m)}(w_{n})
\right)
$$
for $n \leq 2m$.
Note that in the above sum there are no words of ${\cal A}_{1}$
or ${\cal A}_{2}$ at sites greater than $m$. They are only 
occupied by powers of $t$, but then $\widetilde{\Phi}^{(s)}$ sends
them into $1$'s. Therefore, each $j_{k_{i},m_{i}}^{(s)}(w_{i})$ can be
replaced by $j_{k_{i},m_{i}}^{(m)}(w_{i})$ and $\widetilde{\Phi}^{(s)}$ by
$\widetilde{\Phi}^{(m)}$. This ends the proof.
\hfill $\Box$\\
\indent{\par}
{\sc Corollary 4.3.}
{\it $\widetilde{\Phi}^{(m)}\circ j^{(m)}$ converges pointwise to the
conditionally free product od states.}\\[5pt]
{\it Proof.} Obvious.

An extension of the construction presented above
to the case of infinitely many free *-algebras is
very natural. We will show how to do the construction of
$m$-freeness, but we will not repeat the proofs since they
require only minor modifications.

Let ${\cal A}_{l}$, $l\in {\bf N}$ be a family of
unital free *-algebras generated by ${\cal G}_{l}^{+}$. 
Let ${\cal G}_{l}^{-}=\{a^{*}|a\in {\cal G}_{l}^{+}\}$,
${\cal G}_{l}={\cal G}_{l}^{-}\cup {\cal G}_{l}^{+}$.
As before, denote by $\widetilde{{\cal A}}_{l}={\cal
A}_{l}*{\bf C}[t]$ the free product of ${\cal A}_{l}$ and the
algebra of polynomials in one variable $t$. For ach $l\in {\bf
N}$ we identify the units of ${\cal A}_{l}$ and ${\bf C}[t]$
and, by abuse of notation, we denote the unit of each such 
product by ${\bf 1}$. 
As before, extend states $\phi_{l}, \psi_{l}$ on ${\cal 
A}_{l}$ to $\widetilde{\phi}_{l}, \widetilde{\phi}_{l}$ on
$\widetilde{\cal A}_{l}$, $l\in {\bf N}$. 
In the free product $*_{l\in {\bf N}}{\cal A}_{l}$ we
identify the units of ${\cal A}_{l}$, $l\in {\bf N}$.
Abusing notation, we also in this case denote the sequences of
unital free *-algebras, *-homomorphisms and states by 
$j^{(m)}$, ${\cal A}^{(m)}$ and $\widetilde{\Phi}^{(m)}$,
respectively.  \\
\indent{\par}
{\sc Definition 4.4.}
{\it For given $a\in {\cal G}_{l}$, let}
$$
j^{(m)}_{l}(a)=
\sum_{k=1}^{m}
t_{[k,m]}^{\otimes (l-1)}\otimes 
(i_{k,m}(a)-i_{k+1,m}(a))
\otimes 
t_{[k,m]}^{\otimes \infty}
$$
{\it where $t_{[k,m]}^{\otimes \infty}\equiv
(t_{[k,m]})^{\otimes \infty}$ and $i_{m+1,m}(a)=0$,
and define the *-homomorphism}
$$
j^{(m)}:\;\; 
*_{l\in {\bf N}}{\cal A}_{l}\rightarrow 
\bigotimes _{l\in {\bf N}}
\widetilde{\cal A}_{l}^{\otimes m}
$$
{\it as the linear extension of $j^{(m)}({\bf 1})=I_{m}^{\otimes \infty}$ and}
$$
j^{(m)}(w_{1}\ldots w_{n})=j^{(m)}_{k_{1}}(w_{1})\ldots 
j^{(m)}_{k_{n}}(w_{n})
$$
{\it where $w_{1}, \ldots, w_{n}$ are non-empty words in ${\cal A}_{k_{1}},
\ldots, {\cal A}_{k_{n}}$ with $k_{1}, \ldots k_{n} \in {\bf N}$}.

Consider the noncommutative probability space 
($\widetilde{{\cal A}}^{(m)}, \widetilde{\Phi}^{(m)}$), where
$$
\widetilde{{\cal A}}^{(m)}=
\bigotimes_{i\in {\bf N}}
\widetilde{{\cal A}}_{i}^{\otimes m}
$$
and the state is given by
$$
\widetilde{\Phi}^{(m)}=\bigotimes_{i\in {\bf
N}}\widetilde{\Phi}^{(m)}_{i} ,\;\;
\widetilde{\Phi}^{(m)}_{i}=\widetilde{\phi}_{i}\otimes
{\psi}_{i}^{\otimes (m-1)}.
$$
We will also use the $m$-th ``condition'' maps
$$
{\Psi}^{(m)}=\bigotimes_{i\in {\bf
N}}{\Psi}^{(m)}_{i},\;\;\;
\Psi^{(m)}_{i}={\rm id}^{\otimes (m-1)}\otimes \widetilde{\psi}_{i}.
$$
All results of Sections 2-4 are easily generalized to the case
of infinitely many *-algebras. The differences are purely
technical and are omitted. \\
\indent{\par}
{\sc Theorem 4.5.}
{\it $\widetilde{\Phi}^{(m)}\circ j^{(m)}$ agrees with the conditionally
free product on word products $w_{1}\ldots w_{n}$ for $n\leq 2m$, where
$w_{1}, \ldots , w_{n}$ 
$\in {\cal A}_{k_{1}}, \ldots , {\cal A}_{k_{n}}$,
respectively, and $k_{1}\neq k_{2}\neq \ldots \neq k_{n}$,
$k_{1}, \ldots , k_{n}\in {\bf N}$. Thus, $\widetilde{\Phi}^{(m)}\circ
j^{(m)}$ converges pointwise to the conditionally free product
of states.}

If we want to consider an uncountable
number of *-algebras, we can take the continuous tensor product
and proceed in a similar way.\\[5pt]
\begin{center}
{\sc 5. Construction of the Associated *-Bialgebras}
\end{center}
The tensor product constructions are good enough as long as
we only want to study independence of certain variables. 
However, we also would like to
associate a *-bialgebra with each kind of independence.
In the case of the conditionally free independence 
it seems that one should be able to do that for each $m\in {\bf N}$
using the $m$-fold tensor product 
$\widetilde{{\cal A}}^{\otimes m}$. Nevertheless,
it turns out that one needs to take 
the $m$-fold free product $\widetilde{{\cal A}}^{*(m)}$.
The construction of this *-bialgebra is presented below.
First, it is convenient to introduce a free version of the
*-homomorphism $j^{(m)}$.\\
\indent{\par}
{\sc Definition 5.0.}
{\it For given $a\in {\cal G}_{1}, b\in  {\cal G}_{2}$, let
$a_{(i)}$ and $b_{(i)}$, $i\in [m]$, denote different copies of
$a$ and $b$ in $\widetilde{{\cal A}}_{1}^{*(m)}$ and 
$\widetilde{{\cal A}}_{2}^{*(m)}$, respectively, and let $t_{(i)}$, $i\in [m]$,
stand for different copies of $t$ in both products. 
Let}
$$
\widehat{j}_{1}^{(m)}(a)=
\sum_{k=1}^{m}(a_{(k)}-a_{(k+1)})\otimes t_{[k,m]},
$$
$$
\widehat{j}_{2}^{(m)}(b)=\sum_{k=1}^{m}
t_{[k,m]}\otimes (a_{(k)}-a_{(k+1)})
$$
{\it where, by abuse of notation, $t_{[k,m]}=t_{(k)}\ldots t_{(m)}$,
and define the *-homomorphism}
$$
\widehat{j}^{(m)}: {\cal A}_{1}*{\cal A}_{2}\rightarrow 
\widetilde{{\cal A}}_{1}^{* m}
\otimes
\widetilde{{\cal A}}_{2}^{* m}
$$
{\it as the linear extension of $\widehat{j}^{(m)}({\bf 1})={\bf 1}\otimes
{\bf 1}$ and}
$$
\widehat{j}^{(m)}(w_{1}\ldots w_{n})=
\widehat{j}_{k_{1}}^{(m)}(w_{1})\ldots \widehat{j}_{k_{n}}^{(m)}(w_{n}),
$$
{\it where $w_{1}, \ldots , w_{n}$ are non-empty words in
${\cal A}_{k_{1}}, \ldots , {\cal A}_{k_{n}}$,
where $k_{1},\ldots, k_{n}\in \{1,2\}$.}

Let ${\cal A}_{1}={\cal A}_{2}={\cal A}$ in the above definition.
We can associate a cocommutative *-bialgebra with the pair 
(${\cal A},\widehat{j}^{(m)}$). 
Also, let
$\delta:\;{\cal A}\rightarrow {\cal A}*{\cal A}$ be
the *-homomorphism defined by $\delta ({\bf 1})={\bf 1}$, $\delta
(a)=a^{(1)}+a^{(2)}$, where $a^{(1)}, a^{(2)}$ 
are different copies of $a\in {\cal G}$ in ${\cal A}*{\cal A}$. 
Note that $\delta$ maps a given $a$ to the sum of different
copies of $a$. Thus the moments of $\delta(a)$ in the product
state are the moments of the sum of ``independent'', identically
distributed random variables.
\\
\indent{\par}
{\sc Theorem 5.1.}
{\it The *-algebra $\widetilde{{\cal A}}^{*m}$ 
can be equipped with the coproduct }
$$
\Delta^{(m)}({\bf 1})={\bf 1}\otimes {\bf 1},\;\;
\Delta^{(m)}(t_{(i)})=t_{(i)}\otimes t_{(i)},
$$
$$
\Delta^{(m)}(a_{(k)}-a_{(k+1)})=(a_{(k)}-a_{(k+1)})\otimes t_{[k,m]} 
+t_{[k,m]}\otimes (a_{(k)}-a_{(k+1)})
$$
{\it where $k\in [m]$ and it is understood that $a_{(m+1)}=0$,
and the counit}
$$
\epsilon^{(m)} (t_{(k)})= \epsilon^{(m)} ({\bf 1})=1,\;\;
\epsilon^{(m)}(a_{(k)})=0.
$$
{\it Moreover,}
$$
\widehat{j}^{(m)}\circ \delta =\Delta^{(m)} \circ \widehat{i}_{1}
$$
{\it where $\widehat{i}_{1}:\;{\cal A}\rightarrow \widetilde{{\cal A}}^{*(m)}$ is
the canonical *-homomorphic embedding given by
$\widehat{i}_{1}(a)=a_{(1)}$.}\\[5pt]
{\it Proof.}
Note that $t_{(k)}$, $k\in [m]$ (and thus also $t_{[k,m]}$,
$k\in [m]$) are 
group-like and $a_{(k)}-a_{(k+1)}$, $k\in [m]$ are
$t_{[k,m]}$-primitive. Thus it is easy to see that 
$\Delta^{(m)}$ is coassociative.
Verifying that $\epsilon^{(m)}$ is the counit is also immediate.
Therefore 
($\widetilde{{\cal A}}^{*(m)},\Delta^{(m)} , \epsilon^{(m)}$) 
becomes a *-bialgebra. 

Now, $(\widehat{j}^{(m)}\circ \delta)({\bf 1})= 
{\bf 1}\otimes {\bf 1}=(\Delta^{(m)} \circ
\widehat{i}_{1})(1) $. If $a \in {\cal G}$, then
$$
(\widehat{j}^{(m)} \circ \delta)(a)=\widehat{j}_{1}^{(m)}(a)+\widehat{j}_{2}^{(m)}(a)
$$
$$
=
\sum_{k=1}^{m}(a_{(k)}-a_{(k+1)})\otimes t_{[k,m]}
+
\sum_{k=1}^{m}t_{[k,m]}\otimes (a_{(k)}-a_{(k+1)})
$$
$$
=\Delta^{(m)} \left((a_{(1)}-a_{(2)})+
\ldots + (a_{(m-1)}-a_{(m)})+a_{(m)}\right)
=\Delta^{(m)}(a_{(1)})=\Delta^{(m)}\circ \widehat{i}_{1}(a).
$$
This implies that this identity
holds also for arbitrary words in ${\cal A}$ since 
$\widehat{j}^{(m)},\Delta^{(m)}, \delta $ and $\widehat{i}_{1}$ are
*-homomorphisms. 
\hfill
$\Box$

The above theorem shows a relation between
$\widehat{j}^{(m)}\circ \delta$ and the coproduct $\Delta^{(m)}$.
Namely, $\widehat{j}^{(m)} \circ \delta $ equals the coproduct 
$\Delta^{(m)}$ when restricted to the *-subalgebra
$\widehat{i}_{1}({\cal A})$. However, $\Delta^{(m)}$ takes 
$\widehat{i}_{1}({\cal A})$ out of 
$\widehat{i}_{1}({\cal A})\otimes \widehat{i}_{1}({\cal A})$. 
That is why we had to take a bigger *-bialgebra.
Moreover, looking at $\widehat{j}^{(m)}\circ \delta (a)$, we can see that
the $m$-fold free product $\widetilde{{\cal A}}^{*(m)}$ is the right 
choice if we do not assume any additional relations between
different copies of the generators of $\widetilde{{\cal A}}$. 
If we do, we can take the quotient of $\widetilde{{\cal A}}$ modulo
a two-sided ideal providing it is also a coideal. However, many 
relations which appear in the tensor product construction are
not preserved by the coproduct. Thus,
$\Delta^{(m)}$ preserves 
$t_{(i)}t_{(j)}=t_{(j)}t_{(i)}$ for any $i,j$,
and $a_{(i)}t_{(j)}=t_{(j)}a_{(i)}$ for $j<i$, but it does not preserve
$a_{(i)}t_{(j)}=t_{(j)}a_{(i)}$ for $j>i$, or
$a_{(i)}a_{(j)}'=a_{(j)}'a_{(i)}$ for $i\neq j$. 
In other words, the two-sided ideal ${\cal T}$ generated by all
those relations is not a coideal.
This is the reason why we cannot take 
$\widetilde{{\cal A}}^{\otimes m}$ and we have to stick to
$\widetilde{{\cal A}}^{*(m)}$ or its quotient $\widetilde{{\cal
A}}^{*(m)}/{\cal T}_{0}$, where
${\cal T}_{0}$ is the two sided ideal (and a coideal)
generated by $t_{(i)}t_{(j)}=t_{(j)}t_{(i)}$.
All the results can be formulated for either of these two cases
(we choose $\widetilde{\cal A}^{*(m)}$).
\\
\indent{\par}
{\sc Definition 5.2.} 
{\it Let $m,N\in {\bf N}$. For given $a\in {\cal G}_{i}$, $i\in
[N]$ let}
$$
\widehat{j}_{i,N}^{(m)}(a)=\sum_{k=1}^{m}t_{[k,m]}^{\otimes (i-1)}
\otimes (a_{(k)}-a_{(k+1)})\otimes t_{[k,m]}^{\otimes (N-i)},
$$
{\it and define the *-homomorphism}
$$
\widehat{j}_{N}^{(m)}:\;{\cal A}_{1}*\ldots *{\cal A}_{N} \rightarrow 
\bigotimes_{i=1}^{N}
\widetilde{{\cal A}}_{i}^{*(m)}
$$
{\it as the linear extension of $\widehat{j}_{N}^{(m)}({\bf 1})=
{\bf 1}\otimes {\bf 1}$ and}
$$
\widehat{j}_{N}^{(m)}(w_{1}\ldots w_{n})=
\widehat{j}_{k_{1},N}^{(m)}(w_{1})\ldots \widehat{j}_{k_{n},N}^{(m)}(w_{n}),
$$
{\it where $w_{1}, \ldots , w_{n}$ are non-empty words in
${\cal A}_{k_{1}}, \ldots , {\cal A}_{k_{n}}$,
$k_{1},\ldots, k_{n}\in [N]$ and 
$t_{[k,m]}^{\otimes l}\equiv (t_{[k,m]})^{\otimes l}$.}

Then we can express the iterations
of the coproduct in terms of $\widehat{j}_{N}^{(m)}$ in the following way.
\\
\indent{\par}
{\sc Corollary 5.3.}
{\it Let $\delta_{N}$ be the $N$-th iteration of $\delta$, i.e.
$\delta_{N}:\;{\cal A}\rightarrow  
{\cal A} *\ldots *{\cal A}$} ({\it N times}) 
{\it is the *-homomorphism defined by $\delta_{N}({\bf 1})={\bf 1}$,
$\delta_{N}(a)=a^{(1)}+\ldots +a^{(N)}$, where 
$a^{(1)}, \ldots , a^{(N)}$ are different copies of $a$. Then}
$$
\widehat{j}^{(m)}\circ \delta_{N}=\Delta_{N-1}^{(m)}\circ \widehat{i}_{1},
$$
{\it where the $N-1$-th iteration of the coproduct $\Delta^{(m)}$
is obtained from the recursive formula:
$\Delta_{1}^{(m)}=\Delta^{(m)}$,
$\Delta_{k}^{(m)}=({\rm id}\otimes \Delta_{k-1}^{(m)})\circ \Delta^{(m)}$,
$k>1$ . }\\[5pt]
{\it Proof.}
Clearly, 
$$
(\widehat{j}^{(m)}_{N}\circ \delta_{N})({\bf 1})=
{\bf 1}^{\otimes N}=(\Delta_{N-1}^{(m)}\circ \widehat{i}_{1})({\bf 1}).
$$
Now, let $a \in {\cal G}$. Using $\widehat{j}^{(m)}$,
we obtain 
$$
(\widehat{j}^{(m)}_{N}\circ \delta_{N})(a)=
\widehat{j}^{(m)}(a^{(1)}+\ldots +a^{(N)})=
\sum_{i=1}^{N}\widehat{j}_{i,N}^{(m)}(a)
$$
$$
=\sum_{i=1}^{N}\sum_{k=1}^{m}t_{[k,m]}^{\otimes (i-1)}\otimes
(a_{(k)}-a_{(k+1)}) \otimes t_{[k,m]}^{\otimes (N-i)}
$$
$$
=\sum_{k=1}^{m}\Delta_{N-1}^{(m)}(a_{(k)}-a_{(k+1)})=
\Delta_{N-1}^{(m)}(a_{(1)})=(\Delta_{N-1}^{(m)}\circ
\widehat{i}_{1})(a). 
$$
Thus $\widehat{j}_{N}^{(m)}\circ \delta_{N}$ and $\Delta_{N-1}^{(m)} \circ
\widehat{i}_{1}$ agree on the unit and the generators of ${\cal A}$. 
This is enough since $\delta_{N}$, $\widehat{j}_{N}^{(m)}$,
$\Delta_{N-1}^{(m)}$ and $\widehat{i}_{1}$ are *-homomorphisms.
\hfill $\Box$\\
\indent{\par}
{\sc Proposition 5.4.}
{\it Let $\bar{\cal A}={\cal A}*{\bf C}[t,t^{-1}]$. 
The *-bialgebra} 
($\widetilde{{\cal A}}^{*(m)}, \Delta^{(m)},\epsilon^{(m)}$)
{\it can be embedded in the *-Hopf algebra}
($\bar{{\cal A}}^{*(m)}, \bar{\Delta}^{(m)},
\bar{\epsilon}^{(m)}, S^{(m)}$), 
{\it where the coproduct $\bar{\Delta}^{(m)}$ and the counit
$\bar{\epsilon}^{(m)}$ agree with $\Delta^{(m)}$ and $\epsilon^{(m)} $ on
$\widetilde{{\cal A}}^{*(m)}$, respectively, and} 
$$
\bar{\Delta}^{(m)}(t_{(k)}^{-1})=t_{(k)}^{-1}\otimes t_{(k)}^{-1},\;\;\;
\bar{\epsilon}^{(m)}(t_{(k)}^{-1})=1,
$$
{\it with the antipode $S^{(m)}$ defined by}
$$
S^{(m)}({\bf 1})={\bf 1},\;\;
S^{(m)}(t_{(k)})=t_{(k)}^{-1},\;\; 
S^{(m)}(t_{(k)}^{-1})=t_{(k)},
$$
$$
S^{(m)}(a_{(k)}-a_{(k+1)})=-t_{[k,m]}^{-1}(a_{(k)}-a_{(k+1)})t_{[k,m]}^{-1}.
$$
{\it Proof.}
Recalling the definition of a Hopf *-algebra [Kas], we need the
involution and the antipode to satisfy the following conditions: 
(i) * is an antimorphism of real algebras as well as a
morphism of real coalgebras, 
(ii) $S^{(m)}(S^{(m)}(x^{*})^{*})=x$ for all $x\in \bar{{\cal A}}^{*(m)}$. 

The involution * is an antimorphism of real algebras by
definition: 
$(x_{1}\ldots x_{n})^{*}=x_{n}^{*}\ldots x_{1}^{*}$.
To show
that * is a morphism of real coalgebras, one needs
$\bar{\epsilon}^{(m)}$ to be a hermitian functional and the coproduct to
satisfy  
$(*\otimes *)\circ \bar{\Delta}^{(m)}=\bar{\Delta}^{(m)}\circ *$.
The first property follows from the definition of the
counit. Checking the second property for the
generators is immediate. This is enough since
$$
(*\otimes *)\circ \bar{\Delta}^{(m)}(x_{1}\ldots x_{k})=
\bar{\Delta}^{(m)}(x_{k})^{*}\ldots \bar{\Delta}^{(m)}(x_{1})^{*}
$$
$$
=\bar{\Delta}^{(m)}(x_{k}^{*})\ldots \bar{\Delta}^{(m)}(x_{1}^{*})=
\bar{\Delta}^{(m)}(x_{k}^{*}\ldots x_{1}^{*})=
\bar{\Delta}^{(m)}((x_{1}\ldots x_{k})^{*}).
$$
Finally, property (ii) can be easily verified for generators,
from which it follows that it holds for any $x\in \bar{{\cal A}}^{*(m)}$.
\hfill $\Box$

Let us look now at the convolutions of states.
It is known how to define convolutions of states for
*-bialgebras. Namely, if $\Gamma_{2}, \Gamma_{1}$ are two states
on a *-bialgebra (${\cal B}, \Delta, \epsilon$), then the
convolution of $\Gamma_{1}$ and $\Gamma_{2}$ is given by
$$
\Gamma_{1}\star \Gamma_{2}\equiv
(\Gamma_{1}\otimes \Gamma_{2})\circ \Delta.
$$
Thus, in the case of 
$m$-freeness we can express the
convolution of states 
on $\widetilde{{\cal A}}^{*(m)}$ in terms
of the coproduct $\Delta^{(m)}$. 
Let ${\cal T}$ be the two-sided ideal in $\widetilde{{\cal A}}^{*(m)}$
generated by 
$$
t_{(i)}t_{(j)}-t_{(j)}t_{(i)},\;\;
t_{(i)}a_{(j)}-a_{(j)}t_{(i)},\;\;
a_{(i)}a_{(j)}-a_{(j)}a_{(i)},
$$
where $i,j \in [m]$, and $i\neq j$. Then the quotient algebra
$\widetilde{\cal A}^{*(m)}/{\cal T}$ is canonically isomorphic to
$\widetilde{\cal A}^{\otimes m}$. Denote by 
$\eta :\widetilde{\cal A}^{*(m)}\rightarrow $ $\widetilde{\cal
A}^{\otimes m}$ the canonical mapping. Then, for a given state
$\widetilde{\Phi}$ on $\widetilde{\cal A}^{\otimes m}$, let $\widehat{\Phi}=\widetilde{\Phi}
\circ \eta$.  We arrive at the following corollary.\\
\indent{\par}
{\sc Corollary 5.5.}
{\it 
Let $\widetilde{\Phi}^{(m)}_{l}=\widetilde{\phi}_{l}\otimes
\widetilde{\psi}_{l}^{\otimes (m-1)}$, where $\phi_{l},
\psi_{l}$ are states on ${\cal A}$, and let $\widehat{\Phi}^{(m)}_{l}=
\widetilde{\Phi}^{(m)}_{l}\circ \eta$, $l\in [2]$.
Then}
$$
\lim_{m\rightarrow\infty}
(\widehat{\Phi}_{1}^{(m)}\star \widehat{\Phi}_{2}^{(m)})\circ \widehat{i}_{1}=
(\phi_{1},\psi_{1})
\star (\phi_{2}, \psi_{2})
$$
{\it pointwise, 
where $(\phi_{1},\psi_{1})\star (\phi_{2}, \psi_{2})=
*_{i\in \{1,2\}}(\phi_{i},\psi_{i})\circ \delta $
is the conditionally free convolution of states.}\\[5pt]
{\it Proof.}
Let $w$ be a word in ${\cal A}$.
From Theorems 4.0-4.1 and Proposition 5.1 we obtain
$$
\lim_{m\rightarrow \infty}
(\widehat{\Phi}_{1}^{(m)}\star \widehat{\Phi}_{2}^{(m)})
\circ \widehat{i}_{1}(w)
=
\lim_{m\rightarrow \infty}(\widehat{\Phi}_{1}^{(m)}\otimes
\widehat{\Phi}_{2}^{(m)})\circ \Delta^{(m)} \circ \widehat{i}_{1}(w)
$$
$$
=
\lim_{m\rightarrow \infty}
(\widehat{\Phi}_{1}^{(m)}\otimes \widehat{\Phi}_{2}^{(m)})\circ
\widehat{j}^{(m)} \circ \delta (w)
=
*_{i\in \{1,2\}}(\phi_{i},\psi_{i})\circ \delta (w) =
(\phi_{1},\psi_{1})
\star (\phi_{2}, \psi_{2}) (w).
$$
\hfill $\Box$\\
\indent{\par}
{\sc Corollary 5.6.}
{\it Let $\widetilde{\Phi}^{(m)}_{l}=\widetilde{\phi}_{l}\otimes
\widetilde{\psi}_{l}^{\otimes (m-1)}$, where $\phi_{l},
\psi_{l}$ are states on ${\cal A}$, and let $\widehat{\Phi}^{(m)}_{l}=
\widetilde{\Phi}^{(m)}_{l}\circ \eta$, $l\in [N]$.
Then}
$$
\lim_{m\rightarrow \infty}
(\widehat{\Phi}_{1}^{(m)}\star \ldots \star
\widehat{\Phi}_{N}^{(m)}) \circ \widehat{i}_{1}=
(\phi_{1},\psi_{1})\star
\ldots
\star 
(\phi_{N},\psi_{N})
$$
{\it pointwise.}\\[5pt]
{\it Proof.}
In the proof of Corollary 5.5 replace $\delta$ by $\delta_{N}$,
$\Delta^{(m)}$ by $\Delta_{N-1}^{(m)}$ and instead of tensor products and
convolutions of two objects take tensor products and
convolutions, respectively, of $N$ objects.
\hfill $\Box$ \\[5pt]
\begin{center}
{\sc References}
\end{center}
[Av] {\sc D. Avitzour}, {\it Free products of $C^{*}$- algebras}, Trans.
Amer. Math. Soc. {\bf 271} (1982), 423-465.\\[3pt]
[B-L-S] {\sc M.~Bo$\dot{{\rm z}}$ejko, M.~Leinert, R.~Speicher}, 
{\it Convolution and limit theorems for conditionally free
random variables}, Pac. J. Math. 175, No.2 (1996),
357-388.\\[3pt] 
[C-H] {\sc C.D.~Cushen, R.L.~Hudson}, {\it A quantum central limit
theorem}, J. Appl. Prob. {\bf 8}, (1971), 454-469.\\[3pt]
[G-vW] {\sc N.~Giri, W.~von Waldenfels}, {\it An algebraic 
version of the central 
limit theorem}, Z. Wahr. Verw. Gebiete 42 (1978), 129-134.\\[3pt]
[Kas] {\sc Ch.~Kassel}, {\it Quantum groups}, Springer-Verlag, 1995.\\[3pt]
[Len1] {\sc R.~Lenczewski}, {\it On sums of $q$-independent $SU_{q}(2)$ 
quantum variables}, Comm. Math. Phys. 154 (1993), 127-134.\\[3pt]
[Len2] {\sc R.~Lenczewski}, {\it Addition of independent variables in 
quantum groups}, Rev. Math. Phys. 6 (1994), 135-147.\\[3pt] 
[Sch1] {\sc M.~Sch\"{u}rmann}, {\it White Noise on Bialgebras},
Springer-Verlag, Berlin, 1993.\\[3pt]
[Sch2] {\sc M.~Sch\"{u}rmann}, {\it Non-commutative probability on
algebraic structures}, Probability measures on groups and
related structures, Vol. XI (Oberwolfach, 1994), 332-356, World.
Sci. Publishing, River Edge, NJ, 1995.\\[3pt]
[Voi] {\sc D.~Voiculescu}, {\it Symmetries of some reduced free product 
${\cal C}^{*}$-algebras}, in: Operator Algebras and their 
Connections with Topology and Ergodic Theory, Lecture
Notes in Math. 1132, Springer, Berlin, 1985, 556-588.\\[3pt]
[V-D-N] {\sc D.V.~Voiculescu, K.J.~Dykema, A.~Nica}, {\it Free Random
Variables}, CRM Monograph Series, AMS, Providence, 1992.\\[3pt]
[vW] {\sc W.~von Waldenfels}, {\it An approach to the theory of
pressure broadening of spectral lines}, Lecture Note in Math.
{\bf 296} (1973), 19-69.
\end{document}